\documentclass[review]{elsarticle}
\usepackage{xcolor}
\usepackage[latin1]{inputenc}
\usepackage{times}
\usepackage{amssymb,mathrsfs, amsmath}
\usepackage{amsmath}
\allowdisplaybreaks
\usepackage{amsthm}
\usepackage{amssymb}
\usepackage{latexsym}
\usepackage{amsfonts}
\usepackage{epsfig}
\usepackage{hyperref}
\usepackage{graphicx}
\usepackage{graphics}
\usepackage[all]{xy}
\usepackage{tikz-cd}
\usepackage[mathscr]{euscript}
\usepackage{mathrsfs}
\usepackage{stmaryrd}
\usepackage{breqn}
\usepackage{bm}
\usepackage{mathtools}
\usepackage{cancel}

\newtheorem{thm}{Theorem}[section]
\newtheorem{thmm}{Theorem}[subsection]

\newtheorem{prop}{Proposition}[section]

\newtheorem{lem}{Lemma}[section]

\newtheorem{rem}{Remark}[section]

\newtheorem{cor}{Corollary}[section]

\newtheorem{exmpl}{Example}[section]

\theoremstyle{definition}
\newtheorem{defn}{Definition}[section]
\newtheorem{defnn}{Definition}[subsection]

\journal{... }

\DeclareMathOperator{\Hom}{Hom}
\newcommand{\lb }{\llbracket}
\newcommand{\rb}{\rrbracket}
\DeclareMathOperator{\g}{\mathfrak{g}}
\DeclareMathOperator{\h}{\mathfrak{h}}

\begin{document}
\author{Arpan Sharma\corref{myfirstauthor}}
%\address{Sikkim University, Gangtok, Sikkim, 737102, \textsc{India}}
\cortext[myfirstauthor]{This author is supported by INSPIRE, India}
%\cortext[mycorrespondingauthor]{This author is supported by INSPIRE, India}
\ead{rbyadav15@gmail.com, rbyadav01@cus.ac.in}
\author{RB Yadav\corref{mycorrespondingauthor}}
\cortext[mycorrespondingauthor]{Corresponding author}
\address{Sikkim University, Gangtok, Sikkim, 737102, \textsc{India}}
\ead{arpansachin@gmail.com}
\begin{frontmatter}

\title{Crossed homomorphisms on Leibniz superalgebras}

\begin{abstract}

In this article, we give characterisation of  crossed homomorphisms of Leibniz superalgebras as Maurer-Cartan elements of an associated graded Lie algebra. This characterization enables us to develop a cohomology theory for crossed homomorphisms and to apply it to their formal deformation theory. We show that the infinitesimal of a deformation is a 
$1$-cocycle and that linear deformations are classified by such cocycles. Furthermore, we introduce the category of crossed homomorphisms between two fixed Leibniz  superalgebras.
\end{abstract}

\begin{keyword}
\texttt{  Cohomology, Crossed homomorphism, Deformation, Maurer-Cartan element, Leibniz superalgebra}
\MSC[2020]  	16S80, 	16E40,  17B56,  17B10, 17B55,   17A70,   13D10, 13D03  
\end{keyword}
\end{frontmatter}
\section{Introduction}
Leibniz algebra was introduced by Jean-Louis Loday in 1993 \cite{Loday_Leibniz_intro} as a non-commutative generalisation of Lie algebra. Leibniz superalgebras are non-commutative generalisation of Lie superalgebras. Leibniz superalgebra and its cohomologies was introduced by Askar Dzhumadil'dev \cite{MR1888379} in the context of colour Leibniz algebras and was further studied in \cite{Dlui},\cite{MR2102079}. The cohomology and deformation of Leibniz superalgebras were studied in \cite{LeibSup}.

Crossed homomorphisms are also called relative difference operators. Crossed homomorphisms on groups appeared in Whitehead's early work \cite{Whitehead} and were later utilized in the study of non-abelian Galois cohomology \cite{Serre}. The concept was extended to Lie algebras in \cite{Lue} in the context of non-abelian extensions. Crossed homomorphism on Leibniz algebra was seen in \cite{Apurba} whose inverse is said to Rota-Baxter 
operator. More recently, cohomology and deformation of crossed homomorphisms on various algebras have been studied in \cite{Yunhe1},\cite{Apurba},\cite{Wang},\cite{CCCLS},\cite{yadav2026crossed_asso_sup}.

Organisation of the paper is as follows: In section \ref{Leib_sup_represent} we recall some basic definitions and introduce  Leibniz super $\g$-representation. In section \ref{leibniz_sup_crossed_homo}, we define crossed homomorphisms on Leibniz superalgebrs and give some examples. We also construct categotry of crossed homomorphisms on Leibniz superalgebras. In setction \ref{leib_sup_cohomology}, we construct differential graded Lie algebra(dgLa) and characterise crossed homomorphisms as Maurer-Cartan elements of the dgLa and further define its cohomology. In section \ref{leib_sup_deform}, we study linear and formal deformation of crossed homomorphisms on Leibniz superalgebras and establish the corresponding obstruction theory. We derive the deformation equations and show that the first-order term of
every formal deformation is a 1-cocycle in the associated cohomology. Furthermore, we prove that the obstruction to extending a deformation of order $N$ to order $N+1$ is a 2-cocycle.
\section{Leibniz superalgebras and its representations}\label{Leib_sup_represent}
We first recall the definitions of Leibniz superalgebra. Throughout the paper we denote a fixed field by $K$. A linear map $f:V \to W$ is said to  be homogeneous of degree $\alpha$ if $f(V_\beta)\subset W_{\alpha+\beta}$, for all  $\beta\in \mathbb{Z}_2= \{0,1\}$. We write $(-1)^{\textrm{degree}(f)}=(-1)^{\alpha_f}$. Elements of $ V_\beta$ are called homogeneous of degree $\beta.$ We write $(-1)^{a}$ for $(-1)^{\textrm{degree(a)}}$, where $a$ is a homogeneous element of $V$.

\begin{defn} \label{Leibniz_defn} \cite{Dlui}
    A (left) Leibniz superalgebra is a $\mathbb{Z}_2$-graded vector space $\g=\g_0 \oplus \g_1$ equipped with a bilinear map $[\,\,,\,\,]_{\g}:\g \to \g$ satisfying the following conditions:
    \begin{enumerate}
        \item $[a,b]_{\g} \in \g_{\alpha+\beta}$,
        \item $[[a,b]_{\g},c]_{\g}=[a,[b,c]_{\g}]_{\g} -(-1)^{\alpha \beta}[b,[a,c]_{\g}]_{\g}$,
    \end{enumerate}
    for all $a \in \g_\alpha, b \in \g_\beta$ and $c \in \g_\gamma$. If the second condition in the definition \ref{Leibniz_defn} is replaced by $[a,[b,c]_{\g}]_{\g}=[[a,b]_{\g},c]_{\g} -(-1)^{bc}[[a,c]_{\g},b]_{\g}$, then $\g$ is called right Leibniz superalgebra. Throughout this paper, by Leibniz superalgebra we mean it to be left Leibniz superalgebra unless stated otherwise.
\end{defn}

\begin{defn}
    Let $(\g,[-,-]_{\g})$ and $(\h,[-,-]_{\h})$
 be two Leibniz superalgebras. A \textbf{homomorphism} $f: \g \to \h$ is a $K$-linear map such that $f([a,b]_{\g})=[f(a),f(b)]_{\h}$.
 \end{defn}

\begin{defn}\label{representation_Leibniz}\cite{Dlui}
    A representation of a Leibniz superalgebra $(\g,[\,\,,\,\,]_{\g})$ consists of a triple $(V,\rho^L,\rho^R)$ of a $\mathbb{Z}_2$-graded vector space $V$ with two linear maps $\rho^L,\rho^R: \g \to gl(V)$ of degree $0$ satisfying for all homogeneous $a,b \in \g,m \in V$,
    \begin{align}
        \rho^L([a,b]_{\g})(m)&= \rho^L(a)\rho^L(b)(m)-(-1)^{ab}\rho^L(b)\rho^L(a)(m) \label{rep1}\\
        \rho^R(b)\rho^L(a)(m)&=\rho^L(a)\rho^R(b)(m)-(-1)^{am}\rho^R([a,b]_{\g})(m) \label{rep2}\\
        \rho^R(b)\rho^R(a)(m)&=\rho^R([a,b]_{\g})(m)-(-1)^{ma}\rho^L(a)\rho^R(b)(m) \label{rep3}
    \end{align}
From $\eqref{rep2}$ and $\eqref{rep3}$ we get $ \rho^R(b)\rho^L(a)(m)=-(-1)^{am}\rho^R(b)\rho^R(a)(m)$.
In particular, any Leibniz superalgebra $\g$ is a representation of itself with
\[
\rho^L(a)=[a,-]_{\g}, \rho^R(a)=[-,a]_{\g} \mbox{ for } a \in \g.
\]
We call this as regular representation of $
\g$.
\end{defn}
There is another equivalent definition for representation that is given as follows: 
\begin{defn}\cite{MR2286721}
  Let $\g=\g_0\oplus \g_1$ be a Leibniz superalgebra with Leibniz superalgebra bracket denoted by $\pi$. A $\mathbb{Z}_2$-graded vector space $V=V_0\oplus V_1$ over the field K is called a representation of $\g$ if there exist two linear maps (left) $\rho^L:\g \otimes V \to V$ and (right) $\rho^R:V \otimes \g \to V$ satisfying for all homogeneous $a,b \in \g,m \in V$,
  \begin{align}
        \rho^L(\pi(a\otimes b)\otimes m)&= \rho^L(a\otimes \rho^L(b\otimes m)-(-1)^{ab}\rho^L(b\otimes \rho^L(a\otimes m) \label{rep11}\\
        \rho^R( \rho^L(a\otimes m)\otimes b)&=\rho^L(a\otimes \rho^R(m\otimes b))-(-1)^{am}\rho^R(m\otimes \pi(a\otimes b)) \label{rep22}\\
        \rho^R(\rho^R(m\otimes a)\otimes b)&=\rho^R(m\otimes \pi(a\otimes b))-(-1)^{ma}\rho^L(a\otimes \rho^R(m\otimes b)) \label{rep33}
    \end{align}
    \end{defn}
\begin{defn}\label{Leibniz_super_g_representation}
    Let $([\g,[\,\,,\,\,]_\g)$ and $(\h,[\,\,,\,\,]_\h)$ be two Leibniz superalgebras. We say that $\h$ is a Leibniz super $\g$-representation if there are degree $0$ linear maps $\rho^L,\rho^R: \g \to \textrm{End}(\h)$ which makes $(\h,\rho^L,\rho^R)$ into a representation of the Leibniz superalgebra $\g$ and additionally satisfy the following
    \begin{align}
        \rho^L(a)([h,k]_\h) &=[\rho^L(a)(h),k]_\h + (-1)^{ha}[h,\rho^L(a)(k)]_\h, \label{rep4}\\
        [h,\rho^R(a)(k)]_\h &=\rho^R(a)([h,k]_\h) +(-1)^{kh}[k,\rho^R(a)(h)]_\h,\label{rep5}\\
        [h,\rho^L(a)(k)]_\h &=[\rho^R(a)(h),k]_\h +(-1)^{ah}\rho^L(a)([h,k]_\h) \label{rep6}
    \end{align}
    for homogeneous $x \in \g$ and homogeneous $h,k \in \h$.
From $\eqref{rep4}$ and $\eqref{rep6}$ we get $[\rho^L(a)(h),k]_\h=-(-1)^{ah}[\rho^R(a)(h),k]_\h.$ We note that regular representation of any Leibniz superalgebra $\g$ is a Leibniz super $\g$-representation.
\end{defn}

\section{Crossed homomorphisms on Leibniz superalgebras}\label{leibniz_sup_crossed_homo}

\begin{defn} Let $(\g,[-,-]_\g)$ and $(\h,[-,-]_\h)$ be a Leibniz superalgebras and $(\h,\rho^L,\rho^R)$ be Leibniz super $\g$-representation. We define a linear map $D:\g\to \h$ of degree $0$ to be a \textbf{crossed homomorphism} if $$D([x,y]_\g)=\rho^L(x)(Dy)+ \rho^R(y)(Dx)+[Dx, Dy]_\h,\;\forall \; \textrm{ homogeneous }\; x,y \in \g.$$
Crossed homomorphisms are also called as relative difference operators in literature. We call a crossed homomorphism from $\g$ to $\g$ with respect to the adjoint action   a \textbf{difference operator} on $\g$.
We note that if $\rho^L,\rho^R=0,$ then $D$ is simply a Leibniz superalgebra homomorphism.
\end{defn}

%\begin{rem}
  %  Crossed homomorphism on Leibniz algebra \cite{das2022weightedrelative} may %not be a crossed homomorphism on Leibniz superalgebra.
%\end{rem}
%
\begin{exmpl}
    Let $\mathbb{R}^{(1|1)}=\mathbb{R}_0 \oplus\mathbb{R}_1=Span(e_0) \oplus Span(e_1)$, where $e_0=(1,0)$ and $e_1=(0,1)$. $\mathbb{R}^{(1|1)}$ is a Leibniz superalgebra as well as Leibiz algebra with the following bracket:
     \[
     [e_0,e_0]=0,\,\,\,[e_0,e_1]=0,\,\,\, [e_1,e_0]=0, \,\,\, [e_1,e_1]=e_0.
     \]
     With respect to the regular representation of $\mathbb{R}^{(1|1)}$ defined by $\rho^L(x)(y)=[x,y]$ and $\rho^R(x)(y)=[y,x],\,\,\,\textrm{for}\,\,x,y \in \mathbb{R}^{(1|1)}$, suppose $D=\begin{pmatrix}
     a_{11}&a_{12}\\
     a_{21}&a_{22}
     \end{pmatrix}$
     is a crossed homomorphism from $\mathbb{R}^{(1|1)}$ to $\mathbb{R}^{(1|1)}$. Then the following must be satisfied,
     \[
     D([e_i,e_j])=[e_i,De_j]+[De_i,e_j]+[De_i,De_j]
     ,\text{ for all} \; i,j=0,1\]
     $D=\begin{pmatrix}
     -1 & \lambda \\
     0 & -1
     \end{pmatrix}$ for any $\lambda \in \mathbb{R}$ is a solution and it defines a crossed homomorphism on the Leibniz algebra but, for $\lambda \ne 0,$ it is not a crossed homomorphism on the Leibniz superalgebra as its degree is not $0$. Putting $\lambda=0,
     D=\begin{pmatrix}
     -1 & 0 \\
     0 & -1
     \end{pmatrix}
     $ defines a degree $0$ linear map and is a crossed homomorphism on the Leibniz superalgebra.
\end{exmpl}

\begin{exmpl}
    Let us take $\mathfrak{gl}(1|1)_{\mathbb{R}}=M_2(\mathbb{R})$ with the $\mathbb{Z}_2$-grading as
   $\mathfrak{gl}_0=Span\{e_1,e_4\}$ and $\mathfrak{gl}_1=Span\{e_2,e_3\}$ where $e_1=\begin{psmallmatrix}
       1 & 0\\
       0 & 0
   \end{psmallmatrix}, e_2= \begin{psmallmatrix}
       0 & 1 \\
       0 & 0
   \end{psmallmatrix}, e_3= \begin{psmallmatrix}
       0 & 0\\
       1 & 0
   \end{psmallmatrix}, e_4= \begin{psmallmatrix}
       0 & 0\\
       0 & 1
   \end{psmallmatrix}$. With the non-vanishing brackets on $\mathfrak{gl}(1|1)_{\mathbb{R}}$ defined by $[e_1,e_2]=e_2,[e_1,e_3]=e_3,[e_1,e_4]=e_4$ and extending bilinearly makes it a Leibniz superalgebra.
   For any $\lambda \in \mathbb{R}$, $D$ defined by $D(e_1)=-e_1+ \lambda e_4,D(e_2)=-e_2,D(e_3)=-e_3, D(e_4)=-e_4$ gives a one parameter family of crossed homomorphisms on $\mathfrak{gl}(1|1))_{\mathbb{R}}$ with respect to the regular representation.
\end{exmpl}

\begin{lem} \label{induced super representation}
Suppose $(\g,[-,-]_{\g})$ and $(\h,[-,-]_{\h})$ are two Leibniz superalgebras and $(\h,\rho^L,\rho^R)$ is Leibniz super $\g$-representation. Let $D:\g \to \h$ be a crossed homomorphism. Define maps $\rho_D^L,\rho_D^R: \g \to End(\h)$ by
\begin{equation}
\begin{aligned}
\rho_D^L(a)(h)=\rho^L(a)(h) + [Da,h]_{\h}\\
\rho_D^R(a)(h)=\rho^R(a)(h) + [h,Da]_{\h}
\end{aligned}
\end{equation}
for  $a \in \g$ and $h \in \h$. Then, $(\h,\rho_D^L,\rho_D^R)$ is a Leibniz super $\g$-representation.
\end{lem}
\begin{proof}
Clearly, $\rho_D^L$ and $\rho_D^R$ are linear maps of degree $0$. We'll first verify the three representation conditions $\eqref{rep1},\eqref{rep2} \mbox{ and }\eqref{rep3}$. For homogeneous $a,b \in \g$ and $m \in \h$,
\begin{align*}
&\bm{\rho_{D}^{L}(a)\rho_{D}^{L}(b)(m)-(-1)^{ab}\rho_{D}^{L}(b)\rho_{D}^{L}(a)(m)}\\
&=\rho^L(a)(\rho^L(b)(m) + [Db,m]_{\h}) + [Da,\rho^L(b)(m) + [Db,m]_{\h}]_{\h} \\
&\qquad -(-1)^{ab}\rho^L(b)(\rho^L(a)(m) + [Da,m]_{\h})-(-1)^{ab}[Db,\rho^L(a)(m) + [Da,m]_{\h}]_{\h}\\
&= \rho^L(a)\rho^L(b)(m) + \rho^L(a)([Db,m]_{\h}) + [Da,\rho^L(b)(m)]_{\h} + [Da,[Db,m]_{\h}]_{\h}\\
&\qquad -(-1)^{ab}\rho^L(b)\rho^L(a)(m) -(-1)^{ab}\rho^L(b)([Da,m]_{\h})\\ 
&\qquad -(-1)^{ab}[Db,\rho^L(a)(m)]_{\h} -(-1)^{ab}[Db,[Da,m]_{\h}]_{\h}\\
&\mbox{(by \eqref{rep1},\eqref{rep6} and Leibniz identity)}\\
&=\rho^L([a,b]_{\g})(m) + [\rho^L(a)(Db),m]_{\h} +[\rho^R(b)(Da),m]_{\h} + [[Da,Db]_{\h},m]_{\h}\\
&=\rho^L([a,b]_{\g})(m) + [\rho^L(a)(Db) + \rho^R(b)(Da) + [Da,Db]_{\h},m]_{\h}\\
&=\rho^L([a,b]_{\g})(m) + [D([a,b]_{\g}),m]_{\h} \,\,\, \mbox{(D being crossed homomorphism)}\\
&\bm{=\rho_{D}^{L}([a,b]_{\g})(m)}\\ \\
&\bm{\rho_{D}^{L}(a)\rho_{D}^{R}(b)(m) -(-1)^{am}\rho_{D}^{R}([a,b]_{\g})(m)}\\
&=\rho^L(a)\rho^R(b)(m) + \rho^L(a)([m,Db]_{\h}] + [Da,\rho^R(b)(m)]_{\h} + [Da,[m,Db]_{\h}]_{\h}\\
&\qquad -(-1)^{am}\rho^R([a,b]_{\g})(m) -(-1)^{am}[m,D([a,b]_{\g})]_{\h}\\
&\mbox{(by \eqref{rep2} and D being crossed homomorphism)}\\
&=\rho^L(a)\rho^R(b)(m) + \rho^L(a)([m,Db]_{\h}] + [Da,\rho^R(b)(m)]_{\h} + [Da,[m,Db]_{\h}]_{\h}\ 
&\qquad + \rho^R(b)\rho^L(a)(m)-\rho^L(a)\rho^R(b)(m) -(-1)^{am}[m,\rho^L(a)(Db)]_{\h}\\ 
&\qquad -(-1)^{am}[m,\rho^R(b)(Da)]_{\h} -(-1)^{am}[m,[Da,Db]_{\h}]_{\h}\\
&\mbox{(by \eqref{rep4}, \eqref{rep5} and Leibniz identity)}\\
&=\rho^R(b)\rho^L(a)(m) + [\rho^L(a)(m),Db]_{\h} + \rho^R(b)([Da,m]_{\h}) + [[Da,m]_{\h},Db]_{\h}\\
&=\bm{\rho_{D}^{R}(b)\rho_{D}^{L}(a)(m)}\\ \\
&\bm{\rho_{D}^{R}([a,b]_{\g})(m) -(-1)^{am}\rho_{D}^{L}(a)\rho_{D}^{R}(b)(m)}\\
&=\rho^R([a,b]_{\g})(m) + [m,D([a,b]_{\g})]_{\h} -(-1)^{ma}\rho^L(a)\rho^R(b)(m)\\ 
&\qquad -(-1)^{ma}\rho^L(a)([m,Db]_{\h}) -(-1)^{ma}[Da,\rho^R(b)(m)]_{\h} -(-1)^{ma}[Da,[m,Db]_{\h}]_{\h}\\
&\mbox{(by \eqref{rep3} and D being crossed homomorphism)}\\
&=\rho^R(b)\rho^R(a)(m) + (-1)^{am}\rho^L(a)\rho^R(b)(m) + [m,\rho^L(a)(Db)]_{\h} + [m,\rho^R(b)(Da)]_{\h}\\ 
&\qquad + [m,[Da,Db]_{\h}]_{\h} -(-1)^{ma}\rho^L(a)\rho^R(b)(m) -(-1)^{ma}\rho^L(a)([m,Db]_{\h})\\ 
&\qquad-(-1)^{ma}[Da,\rho^R(b)(m)]_{\h} -(-1)^{ma}[Da,[m,Db]_{\h}]_{\h}\\
&\mbox{(by \eqref{rep5},\eqref{rep6} and Leibniz identity)}\\
&=\rho^R(b)\rho^R(a)(m) + \rho^R(b)([m,Da]_{\h}) + [\rho^R(a)(m),Db]_{\h} + [[m,Da]_{\h},Db]_{\h}\\
&\bm{=\rho_{D}^{R}(b)\rho_{D}^{R}(a)(m)} \hspace{4cm}
\end{align*}
Now, we'll verify the additional three Leibniz super $\g$-representation conditions $\eqref{rep4},\eqref{rep5}$ and $\eqref{rep6}$. For homogeneous $a \in \g$ and $h,k \in \h$,
\begin{align*}
&\bm{[\rho_{D}^{L}(a)(h),k]_{\h} + (-1)^{ha}[h,\rho_{D}^{L}(a)(k)]_{\h}}\\
&=[\rho^L(a)(h),k]_{\h} + [[Da,h]_{\h},k]_{\h} + (-1)^{ha}[h,\rho^L(a)(k)]_{\h} + (-1)^{ha}[h,[Da,k]_{\h}]_{\h}\\
&=\rho^L(a)([h,k]_{\h}) + [Da,[h,k]_{\h}]_{\h}\\
&\bm{=\rho_{D}^{L}(a)([h,k]_{\h})} \hspace{4cm}
\end{align*}
By similar computations one can verify the  conditions $\eqref{rep5}$ and $\eqref{rep6}$. Hence $(\h,\rho_D^L,\rho_D^R)$ is a Leibniz super $\g$-representation.
\end{proof}

\begin{lem}\label{semi_direct_Leibniz}
Let $(\h, \rho^L,\rho^R)$ be a Leibniz super $\g$-representation. Then there is a Leibniz superalgebra structure on $\g \oplus \h$ defined by
\begin{equation}
[a+h,b+k]_\ltimes=[a,b]_{\g} + \rho^L(a)(k) + \rho^R(b)(h) + [h,k]_{\h}
\end{equation}
for  $a+h,b+k \in \g\oplus\h$.
\end{lem}
\begin{proof} We have $g\oplus h=(g_0\oplus h_0)\oplus(g_1\oplus h_1).$ This gives a natural $\mathbb{Z}_2$-grading on $g\oplus h$.
Clearly, $[-,-]_\ltimes$ is a degree preserving bilinear map.
For homogeneous $a+h,b+k,c+i \in \g \oplus \h$ such that $|a+h|=\alpha,\, |b+k|=\beta,\, |c+i|=\gamma$,
\begin{align*}
&[a+h,[b+k,c+i]_\ltimes]_\ltimes -(-1)^{\alpha\beta}[b+k,[a+h,c+i]_\ltimes]_\ltimes\\
&=[a,[b,c]_{\g}]_{\g} + \rho^L(a)\rho^L(b)(c) + \rho^L(a)\rho^R(c)(k) + \rho^L(a)([k,i]_{\h})\\
&\qquad + \rho^R([b,c]_{\g})(h) + [h,\rho^L(b)(i)]_{\h} +[h,\rho^R(c)(k)]_{\h} + [h,[k,i]_{\h}]_{\h}\\
&\qquad -(-1)^{\alpha\beta}\Big\{ [b,[a,c]_{\g}]_{\g} + \rho^L(b)\rho^L(a)(i) + \rho^L(b)\rho^R(c)(h) + \rho^L(b)([h,i]_{\h})\\
&\qquad + \rho^R([a,c]_{\g})(k) + [k,\rho^L(a)(i)]_{\h} + [k,\rho^R(c)(h)]_{\h} + [k,[h,i]_{\h}]_{\h} \Big\}\\
&\mbox{(by\eqref{rep1},\eqref{rep2},\eqref{rep3},\eqref{rep4},\eqref{rep5} and \eqref{rep6})}\\
&=[[a,b]_{\g},c]_{\g} + \rho^L([a,b]_{\g})(i) + \rho^R(c)\rho^L(a)(k) + \rho^R(c)\rho^R(b)(h)\\
&\qquad + \rho^R(c)([h,k]_{\h}) + [\rho^L(a)(k),i]_{\h} + [\rho^R(b)(h),i]_{\h} + [[h,k]_{\h},i]_{\h}\\
&=[[a+h,b+k]_\ltimes,c+i]_\ltimes
\end{align*}
\end{proof}
We call this as semi-direct product Leibniz superalgebra and denote it by $\g \ltimes \h$. It follows from lemmas $\ref{induced super representation}$ and $\ref{semi_direct_Leibniz}$, that the direct sum $\g \ltimes \h$ carries another Leibniz superalgebra structure given by
\begin{equation} \label{induced_semi_direct_leibniz}
[a+h,b+k]_{\ltimes_D}=[a,b]_{\g} + \rho_D^L(a)(k) + \rho_D^R(b)(h) + [h,k]_{\h}
\end{equation}
for  $a+h,b+k \in \g\oplus\h$. We denote this semi-direct product Leibniz superalgebra by $\g \ltimes_D \h$. 
\begin{thm}
Suppose $(\g,[-,-]_{\g})$ and $(\h,[-,-]_{\h})$ are two Leibniz superalgebras and $(\h,\rho^L,\rho^R)$ is Leibniz super $\g$-representation. Let $D:\g \to \h$ be a linear map of degree $0$. The following are equivalent;
\begin{enumerate}
    \item $D$ is a crossed homomorphism with repsect to Leibniz super $\g$-representation $(\h,\rho^L,\rho^R)$.
    \item The linear map $\Bar{D}: \g \ltimes_D \h \to \g \ltimes \h$ defined by
    \[
    \Bar{D}(a,h)= (a,Da + h)
    \]
    for homogeneous $a \in \g,h \in \h$ is a Leibniz superalgebra isomorphism.
    \item The map $i_D :\g \to \g \ltimes \h$ defined by 
    \[
    i_D(a)=(a,Da)
    \]
    for homogeneous $a \in \g$ is a Leibniz superalgebra homomorphism.
    \item The graph $Gr(D)=\{(a,Da)| a \in \g\}$ is Leibniz subsuperalgebra of $\g \ltimes \h$.
\end{enumerate}
\end{thm}
\begin{proof}
Let $a+h,b+k \in \g \oplus \h$ be homogeneous elements.\\
\textbf{$(1\Leftrightarrow 2)$} It is easy to see that $\Bar{D}$ is a bijection.
\begin{align*}
& \Bar{D}([a+h,b+k]_{\ltimes_D})\\
& =([a,b]_{\g}, D([a,b]_{\g})+ \rho^L(a)(k) + [Da,k]_{\h} + \rho^R(b)(h) + [h,Db]_{\h} + [h,k]_{\h})\\
&[\Bar{D}(a + h),\Bar{D}(b +k)]_{\ltimes}\\
& =([a,b]_{\g},\rho^L(a)(Db) + \rho^L(a)(k)+ \rho^R(b)(Da)+ \rho^R(b)(h) + [Da+h,Db+k]_{\h})
\end{align*}
We can see that $\Bar{D}([a+h,b+k]_{\ltimes_D})=[\Bar{D}(a + h),\Bar{D}(b +k)]_{\ltimes}$ iff $D$ is a crossed homomorphism.\\
\textbf{$(1\Leftrightarrow 3)$}
\begin{align*}
&i_D([a,b]_{\g})=([a,b]_{\g},D([a,b]_{\g}))\\
&[i_Da,i_Db]_{\ltimes}=[(a,Da),(b,Db)]_{\ltimes}=([a,b]_{\g},\rho^L(a)(Db) + \rho^R(b)(Da) + [Da,Db]_{\h})
\end{align*}
We can see that $i_D([a,b]_{\g})=[i_Da,i_Db]_{\ltimes}$ iff $D$ is a crossed homomorphism.\\
\textbf{ $(1\Leftrightarrow 4)$}
\begin{align*}
[(a,Da),(b,Db)]_\ltimes=([a,b])_{\g},\rho^L(a)(Db)+\rho^R(b)(Da)+[Da,Db]_{\h}
\end{align*}
We can see that $[(a,Da),(b,Db)]_\ltimes \in Gr(D)$ iff $D$ is a crossed homomorphism .
\end{proof}

\subsection{Category of crossed homomorphisms}
\begin{defnn} \label{category_Leibniz}
Let $(\h,\rho^L,\rho^R)$ be a Leibniz super $\g$-representation. Let $D,D'$ be two crossed homomorphisms  from $\g$ to $\h$. A morphism from $D$ to $D'$ is a pair $(\phi_1,\phi_2)$ of two Leibniz superalgebra homomorphisms $\phi_1: \g \to \g$ and $\phi_2: \h \to \h$ of degree $0$ such that for all $a\in \g$ and $h\in \h,$
\begin{enumerate}
    \item $D' \circ \phi_1=\phi_2 \circ D$
    \item $\phi_2(\rho^L(a)(h))=\rho^L(\phi_1(a))(\phi_2(h))$
    \item $\phi_2(\rho^R(a)(h))=\rho^R(\phi_1(a))(\phi_2(h))$
\end{enumerate}
If $\phi_1$ and $\phi_2$ are bijective, then we say that  $D$ and $D'$ are isomorphic.
\end{defnn}
We denote the set of morphisms from $D $ to $D'$ by  $Hom(D,D').$
If $(\phi_1,\phi_2)$ is morphism from $D$ to $D'$ and $(\psi_1, \psi_2)$ is a morphism from $D'$ to $D''$, then  one can easily verify that $(\psi_1\circ \phi_1, \psi_2\circ \phi_2)$ is a morphism from $D$ to $D''.$ For $(\phi_1,\phi_2)\in Hom(D,D')$  and $(\psi_1,\psi_2)\in Hom(D',D'')$, we define
\[
    (\psi_1,\psi_2)\circ (\phi_1,\phi_2)= (\psi_1\circ \phi_1,\psi_2\circ \phi_2).
\]
For every crossed homomorphisms $D$ from $\g$ to $\h$, $(I_{\g},I_{\h})$ is a morphism from $D$ to itself. For any morphism $(\phi_1,\phi_2)$ from $D$ to $D'$ we have $(I_{\g},I_{\h}) \circ (\phi_1,\phi_2)= (\phi_1,\phi_2)\circ (I_{\g},I_{\h})=(\phi_1,\phi_2)$.  Hence we have following theorem: 
\begin{thmm}
   The collection  of all crossed homomorphisms from $\g$ to $\h$  forms a category with Hom set given by Definition \ref{category_Leibniz}. 
\end{thmm}
\section{Cohomology of crossed homomorphisms on Leibniz superalgebras}\label{leib_sup_cohomology}
Let $V=V_0\oplus V_1$ be a   $\mathbb{Z}_2$-graded vector space.
 Consider the permutation group $S_n.$ For any $X=(X_1,\ldots, X_n)$ with $X_i\in V_{x_i}$ and $\sigma\in S_n$, define 
$$ K(\sigma, X)= |\{ (i,j): i<j,\; X_{\sigma(i)}\in V_1, X_{\sigma(j)}\in V_1,  \; \sigma(j) < \sigma(i)\} |$$
and $\epsilon(\sigma, X)=\epsilon(\sigma)(-1)^{K(\sigma,X)},$ 
where $| A |$ denotes cardinality of a set $A, $ $\epsilon(\sigma)$ is the signature of $\sigma.$\\
Let $V$ and $W$ be  $\mathbb{Z}_2$-graded  vector spaces.  For each  $n\in \mathbb{N},$ define   $\mathcal{F}_{n,\alpha}(V,W)$ as  the vector space of all homogeneous $n$-linear mappings  $f:V \underset{n\; times}{\underbrace{\times\cdots\times}}V\to W$ of degree $\alpha.$ Denote the set of all $n$-linear mappings  $f:V \underset{n\; times}{\underbrace{\times\cdots\times}}V\to W$ by  $\mathcal{F}_{n}(V,W)$.  Clearly, we have   $\mathcal{F}_{n}(V,W)=\mathcal{F}_{n,0}(V,W)\oplus \mathcal{F}_{n, 1}(V,W)$. Write  $\mathcal{F}_{0}(V,W)=W$ and  $\mathcal{F}_{-n}(V,W)=0$,  $\forall n\in \mathbb{N}.$  Take $\mathcal{F}(V,W)=\bigoplus_{n\in \mathbb{Z}}\mathcal{F}_{n}(V,W)$. 
We know that $\mathcal{F}_n(V,W)\cong  Hom(V^{\otimes n}, W)$. We write  $\mathbf{E}_n= Hom(V^{\otimes n+1}, V)$ and  $\mathbf{E}=\bigoplus_{n\in \mathbb{Z}}\mathbf{E}_n$. If we denote by $\mathbf{E}_{n,\alpha}=\Hom_{\alpha}(V^{\otimes n +1},V)$ as the vector space of all homogeneous elements of $\mathbf{E}_n$ of  degree $\alpha\in \mathbb{Z}_2$, then $\mathbf{E}_n=\mathbf{E}_{n,0}\oplus \mathbf{E}_{n,1}$ and $\mathbf{E}$ is a $\mathbb{Z}\times\mathbb{Z}_2$-graded vector space. Here, we define  $ f\in \Hom(V^{\otimes n +1},V)$ to have   degree $(n,\alpha_f),$ where $\alpha_f \in \mathbb{Z}_2$ if for any homogeneous $x_1,\ldots, x_{n+1}$ in $V$, $f(x_1,\ldots,x_{n+1})$ is a homogeneous element of $V$ and $\deg (f(x_1,\ldots,x_{n+1}))-\sum_i \deg (x_i)$ is congruent to $\alpha_f$ modulo $2$. Define a product $\circ$ on $\mathbf{E}$ as follows:\\
For every $F \in \mathbf{E}_{n,\alpha_f}$ and $F' \in \mathbf{E}_{n',\alpha_{f'}}$
\[
F \circ F'=\sum_{k=1}^{n+1}(-1)^{(k-1)n'}F \circ_k F'
\]
and $\circ_k$ is defined by
\begin{multline*}
F \circ_kF'(x_1,\ldots,x_{n+n'+1}) \\ =\sum_{\sigma \in S_{(k-1,n')}}\epsilon(\sigma,z_k)(-1)^{\alpha_{f'}(x_{\sigma(1)}+ \cdots + x_{\sigma(k-1)})} F(x_{\sigma(1)}, \ldots , x_{\sigma(k-1)},\\
F'(x_{\sigma(k)},\ldots,x_{\sigma(k+n'-1)},x_{k+n'}),x_{k+n'+1},\ldots,x_{n+n'+1}),
\end{multline*}
where $z_k=(x_1,\ldots,x_{k+n'-1})$,  $x_1,\ldots,x_{n+n'+1} $ are homogeneous  elements in $V$ %$z_k=(x_1,\ldots,x_{k+n'-1})$
and $S_{(k-1,n')}$ consists of permutations $\sigma\in S_{(k-1+n')}$ such that  $\sigma(1)<\cdots<\sigma(k-1)$ and  $\sigma(k)<\cdots<\sigma(k-1+n').$

\begin{thm}\label{MLT3}\cite{LeibSup}
$\mathbf{E}$ is a $\mathbb{Z}\times \mathbb{Z}_2$-graded Lie algebra
  with the bracket  $[-,-]$ defined  by 
  \begin{equation}\label{GNR}
      [F,F']=F\circ F'-(-1)^{nn'+\alpha_{f}\alpha_{f'}}F'\circ F
  \end{equation}
  for $F\in \mathbf{E}_{n,\alpha_{f}}$,  $F'\in \mathbf{E}_{n',\alpha_{f'}}$. 
\end{thm}

Let $W_1$ and $W_2$ be two $\mathbb{Z}_2$ graded vector spaces and elements in $W_1$ will be denoted by $x; y; x_i$ and elements in $W_2$ will be denoted by $u; v; v_i$. Let $g:W_2^{\otimes n} \to W_1$ be a linear map. We can construct a linear map $\hat{g} \in Hom((W_1 \oplus W_2)^{\otimes n},W_1 \oplus W_2)$ by
\[
\hat{g}((x_1,v_1) \otimes \cdots \otimes (x_n,v_n))=(g(v_1 \otimes \cdots \otimes v_n),0)
\]
We have \[(W_1 \oplus W_2)^{\otimes n}=\oplus_{\substack{i_j=1, 2\\j=1,\cdots,n}}   W_{i_1}\otimes \cdots \otimes W_{i_n}.\]
In general, for a given linear map $f: W_{i_1}\otimes \cdots \otimes W_{i_n} \to W_j$, where $i_1, \ldots, i_n, j \in \{1,2\}$, we define a linear map $\hat{f} \in Hom((W_1 \oplus W_2)^{\otimes n},W_1 \oplus W_2)$ by
\begin{equation}\label{extension}
\hat{f} = \begin{cases}
   f & on \,\,\,W_{i_1}\otimes \cdots \otimes W_{i_n}\\
   0 & elsewhere
\end{cases}
\end{equation}
For simplicity of the notation we shall denote $\hat{f}$ by $f$ itself and recognize it from context.

\begin{defn}\label{RAdf30}\cite{LeibSup}
We define  a graded derivation of degree $(\alpha,\beta)$ of a $\mathbb{Z}\times \mathbb{Z}_2$-graded Lie algebra $E$ as  a  linear map $\partial :E\to E$ of degree $(\alpha,\beta)$ such that $\partial ([x,y])=[\partial x,y]+(-1)^{\alpha l +\beta m}[x,\partial y]$, for every homogeneous  $x$ of degree $(l,m)$.
A  differential $\mathbb{Z}\times \mathbb{Z}_2$-graded Lie algebra   $(E, [-,-],\partial)$ consists of  a $\mathbb{Z}\times \mathbb{Z}_2$-graded  Lie algebra   $(E, [-,-])$ and a graded derivation $\partial :E\to E$ of degree $ (1,0)$ such that $\partial^2=0$.   Maurer-Cartan elements of $(E, [-,-],\partial)$ is defined as an element $\omega$ of $E$  of degree $(1,0)$ such that $\partial\omega+\frac{1}{2}[\omega,\omega]=0.$ Every $\mathbb{Z}\times \mathbb{Z}_2$-graded Lie algebra $(E,[-,-])$ can be seen as   differential graded Lie algebra with $\partial=0$.
  Hence Maurer-Cartan elements of the $\mathbb{Z}\times \mathbb{Z}_2$-graded  Lie algebra   $E$ is an element $\omega$ of $E$  of degree $(1,0)$ such that $[\omega,\omega]=0.$
\end{defn} 
We have following theorems: 
\begin{thm}\cite{LeibSup}
    Let $\g$ be $\mathbb{Z}_2$-graded vector space. Then, degree $0$ linear map $\pi:\g \otimes \g \to \g $ is a Leibniz superbracket iff $[\pi,\pi]_B=0$, i.e. $\pi$ is a Maurer-Cartan element of the graded Lie algebra $\bigoplus_{n} Hom(\g^{\otimes n},\g)$ with the bracket $[-,-]_B$ defined in  $\ref{GNR}$.
\end{thm}  
\begin{thm}\cite{LeibSup}
   Let $(\g,[-,-]_\g)$ be a Leibniz superalgebra. Then a  $\mathbb{Z}_2$-graded vector space  $V=V_0\oplus V_1$ is a representation  of  $\g$ with left and right actions $\rho^{L}:\g \otimes V \to V$ and $\rho^R: V \otimes \g \to V$ if and only if  $\Theta=\pi + \rho^L+\rho^R$  ($\pi$ is the Leibniz superbracket $[-,-]_\g$) is a Maurer-Cartan element of the graded Lie algebra $(\bigoplus_{n} Hom((\g \oplus V)^{\otimes n},\g \oplus V)$ with the bracket $[-,-]_B$ defined in $\ref{GNR}$.
\end{thm}
%Let $(\g,[-,-]_{\g})$ and $(\h,[-,-]_{\h})$ be two Leibniz superalgebras and $(\h, \rho^L, \rho^R)$ be a Leibniz super $\g$-representation. We denote $[-,-]_{\g}$ and $[-,-]_{\h}$ by $\pi$ and $\mu$ respectively. We identify the action $\rho^L:\g \to gl(\h)$ with an element in $Hom(\g \otimes \h,\h)$ by $\rho^L(x)(v)=\rho^L(x,v)$. Similarly, $\rho^R:\g \to \h$ is identified with an element in $Hom(\h \otimes \g,\h)$ by $\rho^L(x)(v)=\rho^L(v,x)$.
%\begin{thm}
   % $\Theta=\pi + \rho^L+\rho^R$ is a Maurer-Cartan element of the graded Lie algebra $\bigoplus_{n} Hom((\g \oplus \h)^{\otimes n},\g \oplus \h)$ with the bracket $[-,-]$ defined in $\eqref{MLT3}$.
%\end{thm}
%\noindent
Let $(\g,[-,-]_{\g})$ and $(\h,[-,-]_{\h})$ be two Leibniz superalgebras and $(\h, \rho^L, \rho^R)$ be a Leibniz super $\g$-representation. We denote $[-,-]_{\g}$ and $[-,-]_{\h}$ by $\pi$ and $\mu$ respectively. Since $\pi$ is a Maurer-Cartan element of the graded Lie algebra {\footnotesize $\bigoplus_{n} Hom((\g \oplus \h)^{\otimes n},\g \oplus \h)$,} we can define a differential $\partial_\pi=[\pi,-]$ and the derived bracket on the graded space $\bigoplus_{n} Hom(\g^{\otimes n},\h)$ by;
\[
\llbracket f,g \rrbracket=(-1)^{m-1}[[ \mu,f],g]
\]
for $f \in Hom(\g^{\otimes m},\h)$ and $g \in  Hom(\g^{\otimes n},\h)$.\\
We now consider  $ f\in \Hom(\g^{\otimes n},\h)$ to have   degree $(n,\alpha_f),$ where $\alpha_f \in \mathbb{Z}_2$ if for any homogeneous $x_1,\ldots, x_{n}$ in $\g$, $f(x_1,\ldots,x_{n})$ is a homogeneous element of $\h$ and $\deg f(x_1,\ldots,x_{n})-\sum_i \deg (x_i)$ is congruent to $\alpha_f$ modulo $2$.
We define a linear map $d$ of degree $(1,0)$ on $\bigoplus_{n} Hom(\g^{\otimes n},\h)$ by $d=\partial_{\pi+\rho^L+\rho^R}=[\pi+\rho^L+\rho^R,-]$. More explicitly, for $f \in Hom(\g^{\otimes n},\h)$,
\begin{multline*}
df(x_1,\ldots,x_{n+1})=\rho^R(f(x_1,\ldots,x_{n}),x_{n+1})+\\
(-1)^{n-1}\sum\limits_{i=1}^{n}(-1)^{i+1+x_i(x_1+\ldots + x_{i-1}+\alpha_f)}\rho^L(x_i,f(x_1,\dots,\hat{x}_i,\ldots,x_{n+1})) +\\
(-1)^{n-1}\sum\limits_{k=1}^{n}\sum_{i=1}^{k} (-1)^{i+x_i(x_{i+1} + \cdots + x_k)}\\
f(x_1,\ldots,\hat{x}_i,\ldots,x_k,\pi(x_i,x_{k+1}),x_{k+2},\ldots,x_{n+1})
\end{multline*}

\begin{thm} \label{dgLa}
$\left( \bigoplus_{n} \Hom(\g^{\otimes n},\h),\llbracket -, - \rrbracket,d \right)$ is a differential graded Lie algebra.
\end{thm}
\begin{proof}
Let $f \in \Hom(\g^{\otimes m},\h),g \in \Hom(\g^{\otimes n},\h)$ and $h \in \Hom(\g^{\otimes r},\h)$ be homogeneous.
\begin{align*}
\llbracket f,g \rrbracket &=(-1)^{m-1}[[\mu,f],g]\\
&= (-1)^{m-1} \Big\{ \cancelto{0}{[\mu,[f,g]]} -(-1)^{m-1}[f,[\mu,g]] \Big\}\\
&= - [f,[\mu,g]]\\
&= (-1)^{(m-1)n + \alpha_f\alpha_g}[[\mu,g],f]\\
&= - (-1)^{mn + \alpha_f\alpha_g}\llbracket g,f \rrbracket \hspace{0.5cm} \mbox{(skew symmetry)}\\
\llbracket f,\llbracket g,h \rrbracket \rrbracket &=(-1)^{m+n}[[\mu,f],[[\mu,g],h]]\\
&=(-1)^{m+n}\Big\{[[[\mu,f],[\mu,g]],h]+(-1)^{mn+ \alpha_f\alpha_g}[[\mu,g],[[\mu,f],h]] \Big\}\\
&=(-1)^{m+n}\Big \{\cancelto{0}{[[[[\mu,f],\mu],g],h]} +(-1)^m[[\mu,[[\mu,f],g]],h] \\ 
&\qquad\qquad\qquad\qquad +(-1)^{mn + \alpha_f\alpha_g}[[\mu,g],[[\mu,f],h]]\Big\}\\
&=(-1)^n [[\mu,[[\mu,f],g]],h] +(-1)^{m+n+mn+ \alpha_f\alpha_g}[[\mu,g],[[\mu,f],h]]\\
&=\llbracket\llbracket f,g \rrbracket,h\rrbracket +(-1)^{mn + \alpha_f\alpha_g} \llbracket g,\llbracket f,h \rrbracket \rrbracket \hspace{1cm} \mbox{(Jacobi identity)}
\end{align*}
Clearly, $d$ is a linear map of degree $(1,0)$. By using the fact that $\pi + \rho^L+\rho^R$ is a Maurer-Cartan element of the graded Lie algebra $\Big( \bigoplus_{n} Hom((\g \oplus \h)^{\otimes n},\g \oplus \h),[-,-] \Big)$, it can be easily shown that $d^2=0$. Now, for homogeneous $f \in \Hom(\g^{\otimes m},\h)$ and $g \in \Hom(\g^{\otimes n},\h)$,
\begin{align*}
& d(\llbracket f,g \rrbracket)=(-1)^{m-1}[\pi + \rho^L + \rho^R,[[\mu,f],g]]\\
&=(-1)^{m-1}\Big\{ [[\pi+\rho^L+\rho^R,[\mu,f]],g] + (-1)^m [[\mu,f],[[\pi+\rho^L+\rho^R],g]] \Big\}\\
&=(-1)^{m-1}\Big\{ [[[\pi+\rho^L+\rho^R,\mu],f]-[\mu,[\pi+\rho^L+\rho^R,f]],g] \\
&\qquad\qquad\qquad\qquad + (-1)^m [[\mu,f],[[\pi+\rho^L+\rho^R],g]] \Big\}\\
&=(-1)^m[[\mu,[\pi+\rho^L+\rho^R,f]],g] - [[\mu,f],[[\pi+\rho^L+\rho^R],g]]\\
&=\llbracket[\pi+\rho^L+\rho^R,f],g \rrbracket _ (-1)^{m-1} \llbracket f,[\pi+\rho^L+\rho^R,g] \rrbracket\\
&=\llbracket df,g \rrbracket + (-1)^m\llbracket f,dg \rrbracket
\end{align*}
Hence, $\left( \bigoplus_{n} \Hom(\g^{\otimes n},\h),\llbracket -, - \rrbracket,d \right)$ is a differential graded Lie algebra.
\end{proof}

\begin{thm}
A linear map $D: \g \to \h$ of degree $0$ is a crossed homomorphism iff $D$ is a Maurer-Cartan element of the dgLa $\left( \bigoplus_{n} \Hom(\g^{\otimes n},\h),\llbracket -, - \rrbracket,d \right)$.
\end{thm}
\begin{proof}
Given a linear map $D: \g \to \h$ of degree $0$ is an element in $\bigoplus_{n} \Hom(\g^{\otimes n},\h)$ of degree $(1,0)$. For any $x_1, x_2 \in \g$,
\begin{align*}
&(dD + \frac{1}{2}\llbracket D,D \rrbracket)(x_1,x_2)\\
&\qquad=[\pi+\rho^L + \rho^R,D](x_1,x_2) + \frac{1}{2}\llbracket D,D \rrbracket (x_1,x_2)\\
&\qquad= -D(\pi(x_1,x_2)) + \rho^L(x_1,Dx_2) + \rho^R(Dx_1,x_2) + \mu(Dx_1,Dx_2)
\end{align*}
Hence, $D$ is a crossed homomrphism iff it is a Maurer-Cartan element of the dgLa $\left( \bigoplus_{n} \Hom(\g^{\otimes n},\h),\llbracket -, - \rrbracket,d \right)$.
\end{proof}
A crossed homomorphism $D: \g \to \h$ being Maurer-Cartan element induces a differential $d_D=d+\llbracket D,- \rrbracket$ on the graded Lie algebra $\left( \bigoplus_{n} \Hom(\g^{\otimes n},\h),\llbracket -, - \rrbracket \right)$. Define $C^n(\g,\h)=\Hom(\g^{\otimes n},\h)$ and $C^*(\g,\h)=\bigoplus_{n \geq 1} C^n(\g,\h)$. Thus a crossed homomorphism induces a dgLa $\left( C^*(\g,\h),\llbracket -,- \rrbracket ,d_D \right)$.

\begin{thm}
Let $(\g,\pi)$ and $(\h,\mu)$ be two Leibniz superalgebras and $(\h,\rho^L,\rho^R)$ be a Leibniz super $\g$-representation. Suppose $D:\g \to \h$ is a crossed homomorphism. Then for any linear map $D':\g \to \h$, the sum $D+D'$ is a crossed homomorphism iff $D'$ is a Maurer-Cartan element in the dgLa $\left( C^*(\g,\h),\llbracket -,- \rrbracket ,d_D \right)$.
\end{thm}
\begin{proof}
\begin{align*}
d(D+D') + & \frac{1}{2}\llbracket D+D',D+D' \rrbracket\\
&=dD+dD'+ \frac{1}{2} (\llbracket D,D \rrbracket +\llbracket D,D' \rrbracket + \llbracket D',D \rrbracket + \llbracket D',D' \rrbracket)\\
&=dD' + \llbracket D,D' \rrbracket + \frac{1}{2}\llbracket D',D' \rrbracket \\
&=d_DD' + \frac{1}{2}\llbracket D',D' \rrbracket
\end{align*}
\end{proof}
We call $\left( C^*(\g,\h),d_D \right)$ to be the cochain complex of the crossed homomorphism $D$. We denote the set of $n$-cocyles by $\mathcal{Z}^n_D(\g,\h)=\{f \in C^n(\g,\h) \mid d_D(f)=0 \}$, and the set of $n$-coboundaries by $\mathcal{B
}^n_D(\g,\h)=\{d_D(f) \in C^n(\g,\h) \mid f \in C^{n-1}(\g,\h) \}$ and $n^{th}$ cohomology group by
\[
\mathcal{H}_D^n(\g,\h)= \frac{\mathcal{Z}^n_D(\g,\h)}{\mathcal{B
}^n_D(\g,\h)}
\]
The corresponding cohomology group is called the cohomology of the crossed homomorphism $D$ denoted by $\mathcal{H}^*(\g,\h)$.

\section{Deformation of crossed homomorphisms on Leibniz superalgebras}\label{leib_sup_deform}
    Let $(\g,\pi)$ and $(\h,\mu)$  be  any two Leibniz superalgebras.  Let $(\h,\rho^L,\rho^R)$  be a Leibniz super $\g$-representation and $D:\g\to\h$ be a crossed homomorphism. We write $\g[[t]]=\{\sum_{i=0}^\infty x_it^i: x_i\in \g,\; \forall \; i\}$, $\h[[t]]=\{\sum_{i=0}^\infty y_it^i: y_i\in \h,\; \forall \; i\}$.
    \begin{rem}
    \begin{enumerate}
        \item $\pi$ extends to $\tilde{\pi}:\g[[t]]\otimes \g[[t]]\to \g[[t]]$ given by $$\tilde{\pi}(\sum_{i=0}^\infty x_it^i,\sum_{j=1}^\infty y_jt^j)=\sum_{l=1}^\infty\sum_{i+j=l}\pi(x_i,y_j)t^l$$ and $(\g[[t]], \tilde{\pi})$ is a Leibniz superalgebra. Similarly $\mu$ extends to $\tilde{\mu}:\h[[t]]\otimes \h[[t]]\to \h[[t]]$ and  $(\h[[t]], \tilde{\mu})$ is a Leibniz superalgebra.
        \item $\rho^L,\rho^R $ induce linear maps $\tilde{\rho}^L:\g[[t]] \to gl(\h[[t]])$ given by $\tilde{\rho}^L(x)(\sum\limits_{i=1}^\infty y_it^i)=\sum\limits_{i=1}^\infty\rho^L(x) (y_i)t^i)$ and $\tilde{\rho}^R:\g[[t]] \to gl(\h[[t]])$ given by $\tilde{\rho}^R(x)(\sum\limits_{i=1}^\infty y_it^i)=\sum\limits_{i=1}^\infty\rho^R(x) (y_i)t^i)$ such that $(\h[[t]],\tilde{\rho}^L,\tilde{\rho}^R)$ is a Leibniz super $\g[[t]]$-representation.
        \item For simplicity we denote  $\tilde{\pi},\tilde{\mu},\tilde{\rho}^L$ and $\tilde{\rho}^R$ by $\pi,\mu,\rho^L$ and $\rho^R$ respectively and recognize  them from context.
    \end{enumerate}
\end{rem}
\begin{defn}[Linear Deformation]
    A linear one parameter deformation of crossed homomorphism $D$ is a linear map $D_t:\frac{\g[[t]]}{(t^2)} \to \frac{\h[[t]]}{(t^2)}$,  $D_t=D +tD_1$ where $D_1:\g \to \h$ is a linear map of degree $0$, such that the following condition is satisfied,
\begin{equation}\label{lin_deform_Leib}
   D_t[x,y]_{\g}=\rho^L(x)(D_ty) + \rho^R(y)(D_tx)+[D_tx, D_ty]_{\h}.
\end{equation}
 for every homogeneous $x,y\in \g$.
\end{defn}
Substituting $D_t=D+tD_1$ in $\eqref{lin_deform_Leib}$, we get
\begin{align}\label{RA228}
    D_1[x,y]_{\g} &=\rho^L(x)(D_1y) + \rho^R(y)(D_1x)+[D_1x, Dy]_{\h}+[Dx, D_1y]_{\h}\nonumber\\
    &= D_1[x,y]+dD_1(x,y)+ \lb D,D_1\rb(x,y) \nonumber\\
    &=D_1[x,y]+d_DD_1(x,y)
\end{align}
From equation (\ref{RA228}), we conclude following theorem:
 \begin{thm}
     $D+D_1t:\frac{\g[[t]]}{t^2} \to \frac{\h[[t]]}{t^2}$  is a linear one parameter deformation if and only if $D_1$ is a $1$-cocycle in the cochain complex $(C^*(\g,\h),d_D)$ of $D$.
  \end{thm}
\begin{defn}[Formal Deformation]\label{RA225}
A formal one parameter deformation of $D$ is a linear map $D_t: \g[[t]]\to \h[[t]]$ which consists of formal sum
\[
D_t=\sum_{i=0}^\infty D_it^i=D + D_1t + D_2t^2 + \cdots
\]
where $D_i:\g \to \h$ are linear maps of degree $0$ for every $i\ge 0$ and  $D_0=D$, such that $D_t$ is a crosssd homomorphism, i.e. the following condition is  satisfied,
\begin{equation}\label{RA286}
   D_t[x,y]_{\g}=\rho^L(x)(D_ty) + \rho^R(y)(D_tx)+[D_tx, D_ty]_{\h}
\end{equation}
 for every homogeneous $x,y\in \g$.
\end{defn}
 If $D_1\ne 0,$ then we call it  infinitesimal of the deformation $D_t.$  Equation (\ref{RA286}) is equivalent to following:
\begin{equation}\label{RA227}
   D_n[x,y]_{\g}=\rho^L(x)(D_ny) + \rho^R(y)(D_nx)+\sum_{i+j=n}[D_ix, D_jy]_{\h}, \forall n\ge 0.
\end{equation}
In particular, for $n=1$ we have 
\begin{equation}
    D_1[x,y]_{\g} =\rho^L(x)(D_1y) + \rho^R(y)(D_1x)+[D_1x, Dy]_{\h}+[Dx, D_1y]_{\h} \label{Leib_1}
\end{equation}
From \eqref{RA228} and \eqref{Leib_1}, we conclude that  $d_DD_1=0,$ that is, $D_1$ is a cocycle for the cochain complex $(C^*(\g,\h),d_D)$ of the crossed homomorphism $D$.  Hence we have following lemma:
\begin{lem}
    Infinitesimal of a formal one parameter deformation of crossed homomorphism $D$ is a $1$-cocycle in the cochain complex $(C^*(\g,\h),d_D)$ of $D$.
\end{lem}
Now, we rearrange equation \eqref{RA227} as follows,
\begin{align}\label{Leib_2}
-\sum\limits_{\substack{i+j=n\\i,j \ge 1}}[D_ix,D_jy]_{\h} =&-D_n[x,y]_{\g} + \rho^L(x)(D_ny) + \rho^R(y)(D_nx) \\ \nonumber &+ [D_nx,Dy]_{\h} + [Dx,D_ny]_{\h}, \,\,\, \forall n \ge 0 
\end{align}
By direct computations, we get
\begin{align}\label{Leib_5}
\nonumber&-\frac{1}{2}\sum\limits_{\substack{i+j=n\\i,j \ge 1}}\lb D_i,D_j \rb(x,y)\\ \nonumber
&= -\sum\limits_{\substack{i+j=n\\i,j \ge 1}}[D_ix,D_jy]_{\h}\\ \nonumber
&=-D_n[x,y]_{\g} + \rho^L(x)(D_ny) + \rho^R(y)(D_nx) + [D_nx,Dy]_{\h} + [Dx,D_ny]_{\h}\\ \nonumber
&=d_DD_n(x,y) \,\,\, \forall n \ge 0.\\
\mbox{Hence, }&
-\frac{1}{2}\sum\limits_{\substack{i+j=n\\i,j \ge 1}}\lb D_i,D_j \rb=d_DD_n  \,\,\, \forall n \ge 0.
\end{align}
\begin{defn}
A formal one parameter deformation of crossed homomorphism $D$ of order $N$ is a linear map $D_t: \frac{\g[[t]]}{(t^{N+1})}\to \frac{\h[[t]]}{(t^{N+1})}$ which consists of formal sum
\[
D_t=\sum_{i=0}^{N}D_it^i=D + D_1t + D_2t^2 + \cdots + D_Nt^N
\]
where $D_i:\g \to \h$ are linear maps of degree $0$ for every $1\le i \le n$ and  $D_0=D$, such that $D_t$ is a crosssd homomorphism, i.e. the following condition is  satisfied,
\begin{equation}\label{Leib_4}
   D_t[x,y]_{\g}=\rho^L(x)(D_ty) + \rho^R(y)(D_tx)+[D_tx, D_ty]_{\h}.
\end{equation}
 for every homogeneous $x,y\in \g$.\\
 It follows that equation $\eqref{Leib_5}$ must hold for $n=0,1,\ldots,N$, if $D_t=\sum_{i=0}^{N}D_it^i$ is a deformation of order $N$.
\end{defn}
\begin{defn}
A deformation $D_t=D_t=\sum_{i=0}^{N}D_it^i$ of order $N$ is said to be extensible if there is a linear map of degree $0$ $D_{N+1}: \g \to \h$ such that $\tilde{D}=D_t + D_{N+1}t^{N+1}$ is a deformation of order $N+1$.
\end{defn}
Hence, if a deformation $D_t$ of order $N$ is extensible then deformation equation $\eqref{Leib_5}$ has to be satisfied for $n=N+1$, precisely the following equation,
\begin{equation}\label{Leib_6}
-\frac{1}{2}\sum\limits_{\substack{i+j=N+1\\i,j \ge 1}}\lb D_i,D_j \rb=d_DD_{N+1}
\end{equation}
We observe that the left had side of the above equation $\eqref{Leib_6}$ does not involve $D_{N+1}$; it depends on $D_i$ for $i \le N$ only; we call this as \textbf{obstruction} to extend the deformation $D_t$ and denote it by $Ob_{D_t}$ and we call the cohomology class $[Ob_{D_t}] \in \mathcal{H}^2_{D}(\g,\h)$ as \textbf{obstruction class} to extend the deformation $D_t$.
\begin{prop}
The $2$-cochain $Ob_{D_t}$ is a $2$-cocycle in the cohomology of $D$.
\end{prop}
\begin{proof}
\begin{align*}
        &d_D(-\frac{1}{2}\sum\limits_{\substack{i+j=N+1\\i,j \ge 1}}\lb D_i,D_j \rb)\\
        =&-\frac{1}{2}\sum_{\substack{i+j=N+1\\i,j\ge 1}}\{\llbracket d_DD_i,D_j \rrbracket-\llbracket D_i,d_DD_j\rrbracket\}\\
        =&\frac{1}{4}\sum_{\substack{i+j=N+1\\i,j\ge 1}}\llbracket \sum_{\substack{i_1+i_2=i\\i_1,i_2\ge 1}}\llbracket D_{i_1},D_{i_2}\rrbracket,D_j\rrbracket
        -\frac{1}{4}\sum_{\substack{i+j=N+1\\i,j\ge 1}}\llbracket D_i, \sum_{\substack{j_1+j_2=j\\j_1,j_2\ge 1}}\llbracket D_{j_1},D_{j_2}\rrbracket\rrbracket\\
        =&\frac{1}{4}\sum_{\substack{i_1+i_2+j=N+1\\i_1,i_2,j\ge 1}}\llbracket \llbracket D_{i_1},D_{i_2}\rrbracket,D_j\rrbracket +\frac{1}{4}\sum_{\substack{i+j_1+j_2=N+1\\i,j_1,j_2\ge 1}}\llbracket \llbracket D_{j_1},D_{j_2}\rrbracket,D_i\rrbracket\\
        =&\frac{1}{2}\sum_{\substack{i+j+k=N+1\\i,j,k\ge 1}}\llbracket \llbracket D_{i},D_{j}\rrbracket,D_k\rrbracket\\
        =&0 
\end{align*}
\end{proof}
\begin{thm}
A deformation $D_t$ of a crossed homomorphism $D$ of order 
$N$ extends to a deformation of order $N+1$ if and only if the obstruction class $[Ob_{D_t}] \in \mathcal{H}_{D}^2(\g,\h)$ vanishes.
\end{thm}
\begin{cor}
If $\mathcal{H}_{D}^2(\g,\h)=0$ then deformation of $D$ of any finite order $N$ extends to a deformation of order $N+1$.
\end{cor}
\bibliography{arp.bib}

\end{document}